\documentclass[12pt]{article}

\usepackage{amsmath}
\usepackage{amsfonts}
\usepackage{amsthm}
\usepackage{hyperref}

\begin{document}

\newcommand{\A}{\mbox{${{{\cal A}}}$}}
\newcommand{\Ato}{\mbox{ ${{{\overset{\A}{\longrightarrow}}}}$ }}


\author{Attila Losonczi}
\title{Points accessible in average by rearrangement of sequences II}

\date{\today}

\newtheorem{thm}{\qquad Theorem}[section]
\newtheorem{prp}[thm]{\qquad Proposition}
\newtheorem{lem}[thm]{\qquad Lemma}
\newtheorem{cor}[thm]{\qquad Corollary}
\newtheorem{rem}[thm]{\qquad Remark}
\newtheorem{ex}[thm]{\qquad Example}
\newtheorem{df}[thm]{\qquad Definition}
\newtheorem{prb}{\qquad Problem}

\newtheorem*{thm2}{\qquad Theorem}
\newtheorem*{cor2}{\qquad Corollary}
\newtheorem*{prp2}{\qquad Proposition}

\maketitle

\begin{abstract}

\noindent

We continue investigating the set of limit points of averages of rearrangements of a given sequence. First we generalize results from the previous paper: if a sequence is a composition of two sequences, one of which is bounded and one of which tends to infinity, then we show necessary and sufficient condition for expecting non trivial accessible points. 
A new case will be studied too: the sequences composed of two sequences, one tending to $+\infty$, the other tending to $-\infty$.
Then we start to study accumulation points of the averages of rearranged sequences and prove that if a sequence has 4 accumulation points ($a,b,\pm\infty$) then any closed set in $[a,b]$ can be represented as the set of accumulation points of the averages of a certain rearranged sequence.

\noindent
\footnotetext{\noindent
AMS (2010) Subject Classifications:  40A05, 26E60 \\

}

\end{abstract}

\section{Introduction}
In this paper we are going to continue the investigations started in \cite{lapaar}. For basic definitions, examples, ideas, intentions please consult \cite{lapaar}.

In the first part of the paper our main aim is to generalize results from \cite{lapaar}. In \cite{lapaar} we proved several theorems for unbounded sequences where we assumed that the studied sequence is a composition of two sequences, one of which is constantly 0 and one of which tends to infinity. We got theorems like

\begin{thm2}\cite[Theorem 4.3]{lapaar} Let $(a_n)=(b_n)||(c_n)$ where $b_n\equiv 0,\ c_n\to+\infty$. If 
\[\frac{c_n}{\sum\limits_{i=1}^{n-1}c_i}\to 0\]
then $1$ is accessible in average by rearrangement of $(a_n)$.
\end{thm2}

\begin{thm2}\cite[Theorem 4.4]{lapaar} Let $(a_n)=(b_n)||(c_n)$ where $b_n\equiv 0,\ c_n\to+\infty$ and $(c_n)$ is increasing. If $1$ is accessible in average by rearrangement of $(a_n)$ then
\[\frac{c_n}{\sum\limits_{i=1}^{n-1}c_i}\to 0.\]
\end{thm2}

Now we generalize these theorems in two ways: 

1. We will assume only that the ''lower'' subsequence is bounded from above. 

2. We replace 1 with any point above the upper bound of the ''lower'' subsequence.

\medskip

We are also going to investigate a new case when a sequence does not have finite accumulation point but it has both $\pm\infty$ as accumulation points. We present necessary and sufficient condition for non-trivial points accessible in average by rearrangements. 

\smallskip

Previously we described all behavior of $(a_n)$ hence finally we can enumerate all possible cases regarding $AAR_{(a_n)}$. 

\medskip

In the second part of the paper we turn our attention to accumulation points instead of limit points. We investigate what sets we can get as the accumulation points of the averages of rearrangements of the original sequence.

\smallskip

For more details see subsection \ref{ssmr}. 

\subsection{Basic notions and notations}
For easier readability we copy the basic notations from \cite{lapaar}.

Throughout this paper function $\A()$ will denote the arithmetic mean of any number of variables. We will also use the notation $\A(a_i:1\leq i\leq n)$ for $\A(a_1,\dots,a_n)$.
If $H\subset\mathbb{R}$ is a finite set then $\A(H)$ denotes the arithmetic mean of its distinct points.

Let us use the notation $\bar{\mathbb{R}}=\mathbb{R}\cup\{-\infty,+\infty\}$ and consider $\bar{\mathbb{R}}$ as a 2 point compactification of ${\mathbb{R}}$ i.e. a neighborhood base of $+\infty$ is $\{(c,+\infty]:c\in\mathbb{R}\}$.

\begin{df}Let $(a_n)$ be a sequence. We say that $a_n$ tends to $\alpha\in\bar{\mathbb{R}}$ in average if 
\[\lim\limits_{n\to\infty}\A(a_1,\dots,a_n)=\lim\limits_{n\to\infty}\frac{\sum\limits_{i=1}^na_i}{n}=\alpha.\]
We denote it by $a_n\overset{\A}{\longrightarrow}\alpha$. We also use the expression that $\alpha$ is the limit in average of $(a_n)$.
\end{df}

With this notation if a series $\sum a_n$ is Cesaro summable with sum $c$ then we may say that $s_n\overset{\A}{\longrightarrow}c$ where $s_n=\sum\limits_{i=1}^na_i$.

\begin{df}Let $(a_n)$ be a sequence, $\alpha\in\bar{\mathbb{R}}$. We say that $\alpha$ is accessible in average by rearrangement of $(a_n)$ if there exists a rearrangement of $a_n$ i.e. a bijection $p:\mathbb{N}\to\mathbb{N}$ such that $a_{p(n)}\overset{\A}{\longrightarrow}\alpha$.

The set of all such accessible points will be denoted by $AAR_{(a_n)}$.
\end{df}

\begin{df}If $(a_n),(b_n)$ are two sequences then let $(c_n)=(a_n)||(b_n)$ be the sequence defined by $c_{2n}=b_n,c_{2n-1}=a_n\ (n\in\mathbb{N})$.
\end{df}

\begin{df}If $(a_n),(b_n),(c_n)$ are three sequences then we write $(c_n)=(a_n)\cup (b_n)$ if $(a_n), (b_n)$ are distinct subsequences of $(c_n)$ and they altogether cover all elements in $(c_n)$.
\end{df}

\subsection{Brief summary of the main results}\label{ssmr}

We just enumerate some of the most interesting results to give a taste of the topic.

\begin{thm2}Let $(a_n)=(b_n)||(c_n)$ where $b=\varlimsup b_n\in\mathbb{R},\ c_n\to+\infty$ and $(c_n)$ is increasing. If $c>b$ is accessible in average by rearrangement of $(a_n)$ then
\[\frac{c_n}{\sum\limits_{i=1}^{n-1}c_i}\to 0.\]
\end{thm2}

\begin{thm2}Let $(a_n)=(b_n)||(c_n)$ where $\varlimsup b_n=b\in\mathbb{R},\ c_n\to+\infty$. If 
\[\frac{c_n}{\sum\limits_{i=1}^{n-1}c_i}\to 0\]
and $c>b$ then $c$ is accessible in average by rearrangement of $(a_n)$.
\end{thm2}

\begin{thm2}Let $(a_n)=(b_n)||(c_n),\ b_n\to-\infty,\ c_n\to+\infty$. Then $a\in AAR_{(a_n)},\ a\in\mathbb{R}$ if and only if $\varliminf\frac{b_n}{n}=0$ and $\varliminf\frac{c_n}{n}=0$. \qed 
\end{thm2}

\begin{thm2}Let $(a_n)$ be a sequence bounded from below or above, $(a'_n)$ be one of its rearrangements. Let $c_n=\frac{\sum\limits_{i=1}^na'_i}{n}$. Then the set of accumulation points of $(c_n)$ is $[\varliminf c_n,\varlimsup c_n]$.
\end{thm2}

\begin{thm2}Let $(a_n)$ be a sequence such that it has at least 4 accumulation points: $a,b,-\infty,+\infty\ (a,b\in\mathbb{R},a<b)$. Let $Z\subset[a,b]$ be a closed set. Then there is a rearrangement $(a'_n)$ of $(a_n)$ such that the accumulation points of the sequence $(p_n)$ is exactly $Z$ where $p_n=\A(a'_1,\dots,a'_n)$.
\end{thm2}

\section{More on the basic concept}

We introduce terminology for the most frequently used property of sequences.

\begin{df}Let $(c_n)$ be a sequence such that $c_n\to\infty$. We say that $(c_n)$ is balanced if \[\frac{c_n}{\sum\limits_{i=1}^{n-1}c_i}\to 0.\]
\end{df}

We present an equivalent form.
\begin{prp}$(c_n)$ is balanced if and only if
\[\lim\limits_{n\to\infty}\frac{c_{1}}{c_n}+\frac{c_{2}}{c_n}+\dots+\frac{c_{n-1}}{c_n}=\infty.\tag*{\qed}\]
\end{prp}

Let us use the notation $A^{(c_k)}_n=\frac{c_{1}}{c_n}+\frac{c_{2}}{c_n}+\dots+\frac{c_{n-1}}{c_n}$.

\begin{prp}Removing or adding finitely many elements to a sequence does not change the property of being balanced.\qed
\end{prp}

\begin{prp}\label{pcpkb}Let $(c_n),\ (d_n)$ be balanced sequences, $K\in\mathbb{R}$. Then $(K\cdot c_n),\ (c_n+K),\ (c_n+d_n)$ are all balanced as well.
\end{prp}
\begin{proof}$(K\cdot c_n)$ is balanced:
\[\frac{Kc_n}{\sum\limits_{i=1}^{n-1}Kc_i}=\frac{c_n}{\sum\limits_{i=1}^{n-1}c_i}\to 0.\]

$(c_n+K)$ is balanced:
\[\frac{c_n+K}{\sum\limits_{i=1}^{n-1}(c_i+K)}=\frac{c_n+K}{\sum\limits_{i=1}^{n-1}c_i+(n-1)K}=\frac{\frac{c_n}{\sum\limits_{i=1}^{n-1}c_i}+\frac{K}{\sum\limits_{i=1}^{n-1}c_i}}{1+K\frac{n-1}{\sum\limits_{i=1}^{n-1}c_i}}\to 0\]
because $\sum\limits_{i=1}^{n-1}c_i\to\infty$ moreover $\frac{\sum\limits_{i=1}^{n-1}c_i}{n-1}\to\infty$.

$(c_n+d_n)$ is balanced:
\[\frac{\sum\limits_{i=1}^{n-1}(c_i+d_i)}{c_n+d_n}=\frac{c_n\frac{\sum\limits_{i=1}^{n-1}c_i}{c_n}+d_n\frac{\sum\limits_{i=1}^{n-1}d_i}{d_n}}{c_n+d_n}\to+\infty\]
since this is a weighted average of $\frac{\sum\limits_{i=1}^{n-1}c_i}{c_n}$ and $\frac{\sum\limits_{i=1}^{n-1}d_i}{d_n}$, of which both tend to infinity.
\end{proof}

\begin{prp}Let $(c_n)$ be a balanced sequence. Then $\varlimsup\frac{c_{n-1}}{c_n}=1$.
\end{prp}
\begin{proof}Suppose indirectly that $\varlimsup\frac{c_{n-1}}{c_n}=q<1$. Let $q<p<1$. Then there is $N\in\mathbb{N}$ such that $n>N$ implies that $\frac{c_{n-1}}{c_n}<p$. Take such $n$. If $N+1\leq k\leq n-1$ then 
\[\frac{c_{k}}{c_n}=\frac{c_{k}}{c_{k+1}}\cdot\frac{c_{k+1}}{c_{k+2}}\cdots\frac{c_{n-2}}{c_{n-1}}\cdot\frac{c_{n-1}}{c_n}<p^{n-k}.\]
Then 
\[\frac{c_{1}}{c_n}+\dots+\frac{c_{N}}{c_n}+\frac{c_{N+1}}{c_n}+\dots+\frac{c_{n-1}}{c_n}<\frac{c_{1}}{c_n}+\dots+\frac{c_{N}}{c_n}+p^{n-N-1}+\dots+p^2+p<\]
\[\frac{c_{1}}{c_n}+\dots+\frac{c_{N}}{c_n}+p\cdot\frac{1}{1-p}<+\infty\]
which is a contradiction.
\end{proof}

\begin{ex}$(c_n)$ being balanced does not imply that $\frac{c_{n-1}}{c_n}\to 1$ i.e. $\varliminf\frac{c_{n-1}}{c_n}=1$.
\end{ex}
\begin{proof}Let $(c_n)$ be the following sequence: $1,2,2,4,4,4,8,8,8,8,\dots$ i.e. we have $k+1$ pieces of $2^k$ in $(c_n)$.

We show that $(c_n)$ is balanced. Let $c_{n-1}\not=c_{n}=2^k$. Then 
\[A^{(c_n)}_n=\frac{c_{1}}{c_n}+\dots+\frac{c_{n-1}}{c_n}\geq k\cdot\frac{2^{k-1}}{2^k}=\frac{k}{2}.\]
If $c_{n-1}=c_{n}$ then clearly 
\[A^{(c_n)}_{n-1}=\frac{c_{1}}{c_{n-1}}+\frac{c_{2}}{c_{n-1}}+\dots+\frac{c_{n-2}}{c_{n-1}}<\frac{c_{1}}{c_n}+\frac{c_{2}}{c_n}+\dots+\frac{c_{n-1}}{c_n}=A^{(c_n)}_n\]
which together with the previous statement gives that $(c_n)$ is balanced.

Obviously $\frac{c_{n-1}}{c_n}\not\to 1$.
\end{proof}

\begin{prp}If $c_n\to+\infty$ and $\frac{c_{n-1}}{c_n}\to 1$ then $(c_n)$ is balanced.
\end{prp}
\begin{proof}We can assume that $c_n>0$.

Let $M\in\mathbb{N}$ be given. Choose $0<p<1$ such that $p^{2M}>\frac{1}{2}$. Find $N\in\mathbb{N}$ such that $n>N$ implies that $\frac{c_{n-1}}{c_n}>p$. Let $n>N+2M$. If $N+1\leq k\leq n-1$ then 
\[\frac{c_{k}}{c_n}=\frac{c_{k}}{c_{k+1}}\cdot\frac{c_{k+1}}{c_{k+2}}\cdots\frac{c_{n-2}}{c_{n-1}}\cdot\frac{c_{n-1}}{c_n}>p^{n-k}.\]
Then
\[A^{(c_n)}_n=\frac{c_{1}}{c_n}+\dots+\frac{c_{N}}{c_n}+\frac{c_{N+1}}{c_n}+\dots+\frac{c_{n-1}}{c_n}>\frac{c_{N+1}}{c_n}+\dots+\frac{c_{n-1}}{c_n}>\]
\[\frac{c_{n-2M}}{c_n}+\dots+\frac{c_{n-1}}{c_n}>p^{2M}+p^{2M-1}+\dots+p>2M\cdot\frac{1}{2}=M.\qedhere\]
\end{proof}

\begin{prp}If $\big(A^{(c_k)}_n\big)$ increasing and $\lim\limits_{n\to\infty}A^{(c_k)}_n=+\infty$ then $\frac{c_{n-1}}{c_n}\to 1$.
\end{prp}
\begin{proof}Obvious calculation shows that 
\[A^{(c_k)}_{n+1}=A^{(c_k)}_n\frac{c_{n}}{c_{n+1}}+\frac{c_{n}}{c_{n+1}}=(A^{(c_k)}_n+1)\frac{c_{n}}{c_{n+1}}.\]
Then $A^{(c_k)}_{n}\leq A^{(c_k)}_{n+1}$ gives that
\[A^{(c_k)}_{n}\leq (A^{(c_k)}_n+1)\frac{c_{n}}{c_{n+1}}\]
\[A^{(c_k)}_n\big(\frac{c_{n+1}}{c_{n}}-1\big)\leq 1\]
\[A^{(c_k)}_n\leq \frac{1}{\frac{c_{n+1}}{c_{n}}-1}.\]
Which implies that 
\[\frac{1}{\frac{c_{n+1}}{c_{n}}-1}\to+\infty\] 
hence $\frac{c_{n+1}}{c_{n}}-1\to 0$ and $\frac{c_{n+1}}{c_{n}}\to 1$ and finally $\frac{c_{n-1}}{c_{n}}\to 1$.
\end{proof}

\begin{ex}$\lim\limits_{n\to\infty}A^{(c_k)}_n=+\infty$ and $\frac{c_{n-1}}{c_n}\to 1$ do not imply that $\big(A^{(c_k)}_n\big)$ is increasing.
\end{ex}
\begin{proof}Let $(c_n)$ be the following sequence: $1,1,2,2,2,3,3,3,3,...$ where there are $n+1$ pieces from $n$ in the sequence.

Obviously $c_n\to+\infty,\ (c_n)$ is increasing and $\frac{c_{n-1}}{c_n}\to 1$.

We show that $\lim\limits_{n\to\infty}A^{(c_k)}_n=+\infty$. If $c_n=c_{n+1}$ then clearly $A^{(c_k)}_n<A^{(c_k)}_{n+1}$. If $c_n\not=c_{n+1},\ c_{n+1}=k+1$ then $A^{(c_k)}_{n+1}>(k+1)\cdot\frac{k}{k+1}=k$. These two statements gives the claim.

We finally show that $\big(A^{(c_k)}_n\big)$ is not increasing. Let $c_n=k,\ c_{n+1}=k+1$. Then $A^{(c_k)}_{n+1}<A^{(c_k)}_n$ because
\[A^{(c_k)}_n=2\cdot\frac{1}{k}+3\cdot\frac{2}{k}+\dots+k\cdot\frac{k}{k},\]
\[A^{(c_k)}_{n+1}=2\cdot\frac{1}{k+1}+3\cdot\frac{2}{k+1}+\dots+(k+1)\cdot\frac{k}{k+1},\]
and all terms in the latter are smaller than in the former, except the last terms which are equal.
\end{proof}

\begin{ex}$(c_n)$ being balanced does not imply that $(c^2_n)$ is balanced as well.
\end{ex}
\begin{proof}Let $(c_n)$ be the following sequence: 4 pieces of 1 then 9 pieces of 2 then 16 pieces of 6 and so on; generally there are $(n+1)^2$ pieces of $n!$ in the sequence.

We show that $(c_n)$ is balanced. If $c_n=c_{n+1}$ then clearly $A^{(c_k)}_n<A^{(c_k)}_{n+1}$. If $c_n\not=c_{n+1},\ c_{n+1}=k!$ then $A^{(c_k)}_{n+1}>k^2\cdot\frac{(k-1)!}{k!}=k$. These two statements gives the claim.

We show that $(c^2_n)$ is not balanced. Let $c_n\not=c_{n+1},\ c_{n+1}=k!$. Then 
\[A^{(c_k)}_{n+1}=\sum\limits_{l=1}^{k-1}(l+1)^2\cdot\Big(\frac{l!}{k!}\Big)^2=\sum\limits_{l=1}^{k-1}\Big(\frac{(l+1)!}{k!}\Big)^2=\]
\[1+\sum\limits_{l=1}^{k-2}\Big(\frac{(l+1)!}{k!}\Big)^2<1+\sum\limits_{l=1}^{k-2}\frac{1}{k^2}=1+(k-2)\frac{1}{k^2}<1+\frac{1}{k}\leq 2\]
showing that $(c^2_n)$ is not balanced.
\end{proof}

\section{Generalizing previous results}

First let us present a generalized form of \cite[Theorem 3.1]{lapaar}.

\begin{thm}\label{t2t1}Let $(a_n)$ be a sequence and $a\leq b\ (a,b\in\mathbb{R})$ are two of its accumulation points. Then $[a,b]\subset AAR_{(a_n)}$.
\end{thm}
\begin{proof}The proof can copy the proof of \cite[Theorem 3.1]{lapaar} remarking that we have not used that $(d_n)$ is bounded.
\end{proof}

\begin{thm}\label{t1aariff2b}Let $(a_n)=(b_n)||(c_n)$ where $b=\varlimsup b_n\in\mathbb{R},\ c_n\to+\infty$ and $(c_n)$ is increasing. If $c>b$ is accessible in average by rearrangement of $(a_n)$ then
\[\frac{c_n}{\sum\limits_{i=1}^{n-1}c_i}\to 0.\]
\end{thm}
\begin{proof}We can assume that $c>0$ (i.e. $b\geq 0$) because otherwise consider sequences $(b_n-b),(c_n-b)$. If the statement is true for these sequences then by \ref{pcpkb} it is true for the original sequences as well.

Let $(d_n)$ be a rearrangement such that $d_n\overset{\A}{\longrightarrow}c$. This rearrangement defines a rearrangement of $(c_n)$, namely take the elements from $(c_n)$ exactly in the same order as they come in $(d_n)$. Let us denote that rearranged sequence with $(c'_n)$ and $c'_n=d_{m_n}$.

Let $\epsilon>0$. Then there is $N\in\mathbb{N}$ such that $m\geq N$ implies that 
\[c-\epsilon<\frac{\sum\limits_{i=1}^{m}d_i}{m}<c+\epsilon.\]
Let $n$ be chosen such that $m_{n-1}>N$. Set $s_n=\sum\limits_{i=1}^{n-1}c'_i$ 
and \[w_{n-1}=\sum\limits_{i=1}^{m_{n-1}}d_i-\sum\limits_{i=1}^{n-1}c'_i\]
i.e. $w_{n-1}$ is the sum of elements form $(b_i)$ which are among $d_1,\dots,d_{m_{n-1}}$.
We get that
\begin{equation}\label{eq5a2}
c-\epsilon<\frac{\sum\limits_{i=1}^{m_n-1}d_i}{m_n-1}=\frac{s_n+w_n}{m_n-1}<c+\epsilon,
\end{equation}
\begin{equation}\label{eq5b2}
c-\epsilon<\frac{\sum\limits_{i=1}^{m_n}d_i}{m_n}=\frac{s_n+w_n+c'_n}{m_n}<c+\epsilon.
\end{equation}
From (\ref{eq5b2}) we get that
\[c-\epsilon<\frac{1+\frac{c'_n}{s_n+w_n}}{\frac{m_n}{s_n+w_n}}<c+\epsilon.\]
By multiplying with the denominator and using (\ref{eq5a2}) we get that
\[\frac{c-\epsilon}{c+\epsilon}\cdot\frac{m_n}{m_n-1}<(c-\epsilon)\frac{m_n-1}{s_n+w_n}\cdot\frac{m_n}{m_n-1}<1+\frac{c'_n}{s_n+w_n}<\]
\[<(c+\epsilon)\frac{m_n-1}{s_n+w_n}\cdot\frac{m_n}{m_n-1}<\frac{c+\epsilon}{c-\epsilon}\cdot\frac{m_n}{m_n-1}\]
and clearly both sides tend to $1$ when $\epsilon\to 0$.

It gives that $\frac{c'_n}{s_n+w_n}\to 0$.

We show that it implies that $\frac{c'_n}{s_n}\to 0$. Let the number of terms in $w_n$ be $k_n=m_n-n$, the number of terms in $s_n$ be $l_n=n-1$. Clearly \[\frac{c'_n}{s_n+w_n}=\frac{\frac{c'_n}{s_n}}{1+\frac{w_n}{s_n}}.\]
Therefore it is enough to prove that $\frac{w_n}{s_n}$ is bounded. Assume the contrary and assume first that $\varlimsup\frac{w_n}{s_n}=+\infty$. Let $N\in\mathbb{N}$ such that $n>N$ implies that 
\[\frac{s_n+w_n}{l_n+k_n}>\frac{c+b}{2}\text{ and }\frac{w_n}{k_n}<\frac{c+2b}{3}.\]
Let $K\in\mathbb{R}$ such that \[\frac{K+1}{K}\cdot\frac{c+2b}{3}<\frac{c+b}{2}.\] Now choose $n>N$ such that $\frac{w_n}{s_n}>K$. Then \[w_n>Ks_n\]
\[(K+1)w_n>K(s_n+w_n)\]
\[\frac{c+b}{2}>\frac{K+1}{K}\cdot\frac{c+2b}{3}>\frac{K+1}{K}\cdot\frac{w_n}{k_n}>\frac{l_n+k_n}{k_n}\cdot\frac{s_n+w_n}{l_n+k_n}>\frac{s_n+w_n}{l_n+k_n}>\frac{c+b}{2}\]
which is a contradiction.

Assume now that $\varliminf\frac{w_n}{s_n}=-\infty$. Let $N\in\mathbb{N}$ such that $n>N$ implies that 
\[\frac{s_n+w_n}{l_n+k_n}>\frac{c}{2}>0\text{ and }\frac{s_n}{l_n}>0.\]
For $K=-2$ choose $n>N$ such that 
\[\frac{w_n}{s_n}<-2\]
\[w_n<-2s_n\]
\[\frac{s_n+w_n}{l_n+k_n}<\frac{s_n-2s_n}{l_n+k_n}=-\frac{l_n}{l_n+k_n}\frac{s_n}{l_n}<0\]
-- a contradiction again.

Finally we got that $\frac{c'_n}{s_n}\to 0$. Now \cite[Corollary 4.2]{lapaar} yields the statement.
\end{proof}

\begin{thm}\label{t1aariff2}Let $(a_n)=(b_n)||(c_n)$ where $\varlimsup b_n=b\in\mathbb{R},\ c_n\to+\infty$. If 
\[\frac{c_n}{\sum\limits_{i=1}^{n-1}c_i}\to 0\]
and $c>b$ then $c$ is accessible in average by rearrangement of $(a_n)$.
\end{thm}
\begin{proof}Set $v=\frac{1}{c-b}$. Set $s_n=\sum\limits_{i=1}^{n}c_i$. Note that $s_n\to+\infty$ moreover $\frac{s_n}{n}\to+\infty$. With that notation the assumption gets the form: $\frac{c_n}{s_{n-1}}\to 0$. 

Take a subsequence $(b'_n)$ from $(b_n)$ such that $b'_n\to b$. 

We can assume that $c_n$ is increasing (by \cite[Corollary 4.2]{lapaar}) and $c_n>\max\{1,2(c-b)\}$. Then let $m_n=\Big\lfloor v\cdot\sum\limits_{i=1}^{n}c_i \Big\rfloor$. The previous assumption on $c_n$ gives that $(m_n)$ is a strictly increasing sequence of integers.

If $k\in\mathbb{N}$ then let 
\[d_k=\begin{cases}
c_n&\text{if }k=m_n\\
b'_{k-n}&\text{if }k\in(m_n,m_{n+1})
\end{cases}\]
i.e. on the $m_n$ position put $c_n$, on the other positions put elements from $(b'_n)$.  

We show that $d_k\overset{\A}{\longrightarrow}c$. Let $m_{n-1}\leq k<m_n$.
Obviously
\[\A(d_1,\dots,d_k)=\frac{\sum\limits_{i=1}^{k-(n-1)}b'_i+\sum\limits_{i=1}^{n-1}c_i}{k}=\frac{\sum\limits_{i=1}^{k-(n-1)}b'_i}{k}+\frac{\sum\limits_{i=1}^{n-1}c_i}{k}.\]
Let us investigate the two terms. For the first one we get that
\[\frac{\sum\limits_{i=1}^{k-(n-1)}b'_i}{k}=\frac{\sum\limits_{i=1}^{k-(n-1)}b'_i}{k-(n-1)}\cdot\frac{k-(n-1)}{k}.\]
Clearly \[\frac{\sum\limits_{i=1}^{k-(n-1)}b'_i}{k-(n-1)}\to b\]because $b'_i\to b$. Regarding the second factor we have
\[\frac{n-1}{k}\leq\frac{n-1}{m_{n-1}}\leq\frac{n-1}{vs_{n-1}-1}=\frac{1}{v\frac{s_{n-1}}{n-1}-\frac{1}{n-1}}\to 0\]
Hence \[\frac{\sum\limits_{i=1}^{k-(n-1)}b'_i}{k}\to b.\]
Now let us analyze the second term.
\[\frac{1}{v+v\frac{c_n}{s_{n-1}}+\frac{1}{s_{n-1}}}=\frac{s_{n-1}}{v(s_{n-1}+c_n)+1}\leq\frac{s_{n-1}}{m_n}\leq\frac{\sum\limits_{i=1}^{n-1}c_i}{k}\leq\]
\[\leq\frac{s_{n-1}}{m_{n-1}}=\frac{s_{n-1}}{vs_{n-1}-1}=\frac{1}{v+\frac{1}{s_{n-1}}}\]
Evidently both sides tends to $\frac{1}{v}=c-b$ using that $\frac{c_n}{s_{n-1}}\to 0$ by assumption.
Therefore \[\frac{\sum\limits_{i=1}^{n-1}c_i}{k}\to c-b\]
which gives that $\lim\limits_{k\to \infty}\A(d_1,\dots,d_k)=c$

By \cite[Corollary 2.4]{lapaar} we can add the remaining elements from $(a_n)$ (that are not in $(b'_n)||(c_n)$) to $(d_n)$ such that the limit in average does not change.
\end{proof}

\section{Sequences composed of two sequences, one tending to $+\infty$, the other tending to $-\infty$}

We are going to investigate the case when a sequence does not have finite accumulation point but it has both $\pm\infty$ as accumulation points. 

\begin{thm}Let $(a_n)=(b_n)||(c_n),\ b_n\to-\infty,\ c_n\to+\infty$. If there is $a\in AAR_{(a_n)},\ a\in\mathbb{R}$ then $\varliminf\frac{b_n}{n}=0$ and $\varliminf\frac{c_n}{n}=0$.
\end{thm}
\begin{proof}Suppose indirectly that $\varliminf\frac{c_n}{n}=2p>0$ ($(b_n)$ can be handled similarly). Then there is $N\in\mathbb{N}$ such that $n>N$ implies that $\frac{c_n}{n}>p$.

Let $\frac{p}{9}>\epsilon>0$. If $a\in AAR_{(a_n)}$ then there is a rearrangement of $(a_n)$ say $(a_{n_k})$ such that $a_{n_k}\overset{\A}{\longrightarrow}a$. Which yields that there is $M\in\mathbb{N}$ such that $v>M$ implies that 
\[\left\lvert\frac{\sum\limits_{k=1}^va_{n_k}}{v}-a\right\lvert<\epsilon.\]
Let $M$ be chosen such that $\{c_n:n\leq N\}\subset\{a_{n_k}:k\leq M\}$ moreover $\frac{a+1}{M+1}<\epsilon$. Find a $c_n$ such that $n>N$ and $c_n=a_{n_{v+1}}$ and $v>M$ and $v+1<2n+1$. It can be done because let 
\[l=\max\big\{m:(\exists k\leq M\text{ such that } b_m=a_{n_k})\text{ or }(\exists k\leq M\text{ such that } c_m=a_{n_k})\big\}\] 
and let $v+1=\min\{w:a_{n_w}=c_m,\ m>l\}$.

Let us estimate how much the new element $a_{n_{v+1}}$ modifies the average.
\[\left\lvert\frac{\sum\limits_{k=1}^va_{n_k}}{v}-\frac{\sum\limits_{k=1}^{v+1}a_{n_k}}{v+1}\right\lvert=\left\lvert\frac{\sum\limits_{k=1}^va_{n_k}}{v(v+1)}-\frac{c_n}{v+1}.\right\lvert\]
Clearly \[\left\lvert\frac{\sum\limits_{k=1}^va_{n_k}}{v(v+1)}\right\lvert<\frac{a+\epsilon}{v+1}<\frac{a+1}{M+1}<\epsilon\]
and 
\[\frac{c_n}{v+1}>\frac{c_n}{2n+1}>\frac{\frac{c_n}{n}}{2+\frac{1}{n}}>\frac{p}{2+\frac{1}{n}}\geq\frac{p}{3}>3\epsilon\]
which gives that
\[\left\lvert\frac{\sum\limits_{k=1}^va_{n_k}}{v}-\frac{\sum\limits_{k=1}^{v+1}a_{n_k}}{v+1}\right\lvert>2\epsilon\]
which is a contradiction.
\end{proof}

\begin{thm}Let $(a_n)=(b_n)||(c_n),\ b_n\to-\infty,\ c_n\to+\infty$. If $\varliminf\frac{b_n}{n}=0$ and $\varliminf\frac{c_n}{n}=0$ then $\mathbb{R}\subset AAR_{(a_n)}$.
\end{thm}
\begin{proof}If enough to prove that $0\in  AAR_{(a_n)}$ because $(b_n-a)||(c_n-a)$ satisfies the conditions of the theorem and if $a_n-a\overset{\A}{\longrightarrow}0$ then $a_n\overset{\A}{\longrightarrow}a$.

Let $\epsilon>0$.

\smallskip

First assume that $\frac{b_n}{n}\to 0$ and $\frac{c_n}{n}\to 0$. We define a rearrangement of $(a_n)$. First take elements from $(c_n)$ in the order as they are in $(c_n)$ such that the average would be $\geq 0$. It can be done as $c_n\to+\infty$. Then add elements from $(b_n)$ in the order as they are in $(b_n)$ such that the average of all selected elements (the previous ones and the just added ones) would be $\leq 0$. It can be done as $b_n\to-\infty$. Then add not used elements from $(c_n)$ in the order as they are in $(c_n)$ such that the average of all selected elements (the previous ones and the just added ones) would be $\geq 0$. And so on. Let us denote the constructed sequence by $(a_{n_k})$.

\smallskip

We show that the rearranged sequence tends to $0$ in average. Let $N\in\mathbb{N}$ such that $n>N$ implies that $\frac{|b_n|}{n}<\frac{\epsilon}{2},\ \frac{c_n}{n}<\frac{\epsilon}{2}$. 

Choose $V\in\mathbb{N}$ such that $\{b_n,c_n:n\leq N\}\subset\{a_{n_k}:k\leq V\}$. Let 
\[s=\frac{\sum\limits_{k=1}^Va_{n_k}}{V}.\]
Suppose that $s\geq 0$ (the other case is similar).

First we show that there is $v\geq V$ such that \[\frac{\sum\limits_{k=1}^va_{n_k}}{v}\in(-\epsilon,\epsilon).\]
If $s<\epsilon$ then we are done. If not then suppose $\frac{\sum\limits_{k=1}^va_{n_k}}{v}\geq \epsilon$. Then clearly 
\[\frac{\sum\limits_{k=1}^{v+1}a_{n_k}}{v+1}=\frac{\sum\limits_{k=1}^va_{n_k}}{v}\cdot\frac{v}{v+1}+\frac{b_n}{v+1}\]
for $a_{n_{v+1}}=b_n$. But 
\[\frac{\sum\limits_{k=1}^va_{n_k}}{v}\cdot\frac{v}{v+1}>0\text{ and }\frac{|b_n|}{v+1}<\frac{|b_n|}{n}<\frac{\epsilon}{2}\]
hence $\frac{\sum\limits_{k=1}^{v+1}a_{n_k}}{v+1}>-\frac{\epsilon}{2}$ i.e. sooner or later we will step into $(-\epsilon,\epsilon)$ referring also to the definition of $(a_{n_k})$.

\smallskip

Now suppose that $\frac{\sum\limits_{k=1}^va_{n_k}}{v}\in(-\epsilon,\epsilon)$ and let us examine how much the average changes if we add a new element. If $v>V,\ a_{n_{k+1}}=c_n$ ($b_n$ is similar) then
\[\left\lvert\frac{\sum\limits_{k=1}^va_{n_k}}{v}-\frac{\sum\limits_{k=1}^{v+1}a_{n_k}}{v+1}\right\lvert\leq\left\lvert\frac{\sum\limits_{k=1}^va_{n_k}}{v(v+1)}\right\lvert+\left\lvert\frac{c_n}{v+1}\right\lvert\leq\frac{1}{v+1}\cdot\left\lvert\frac{\sum\limits_{k=1}^va_{n_k}}{v}\right\lvert+\left\lvert\frac{c_n}{n}\right\lvert<\frac{\epsilon}{2}+\frac{\epsilon}{2}=\epsilon.\]
Which means that by the definition of $(a_{n_k})$ we will always remain in the interval $(-\epsilon,\epsilon)$. I.e. we showed that $a_{n_k}\overset{\A}{\longrightarrow}0$.

\smallskip

Let us turn to the general case and let $(b_{n_k}),\ (c_{m_k})$ subsequences of $(b_n),\ (c_n)$ respectively such that $b_{n_k}\to 0,\ c_{m_k}\to 0$. Apply the previous assertion for sequence $(d_n)=(b_{n_k})||(c_{m_k})$ which gives that $\mathbb{R}\subset AAR_{(d_n)}$. Now \cite[Corollary 2.4]{lapaar} gives the statement.
\end{proof}

We can summarize the previous two theorems in the following way.

\begin{thm}Let $(a_n)=(b_n)||(c_n),\ b_n\to-\infty,\ c_n\to+\infty$. Then $a\in AAR_{(a_n)},\ a\in\mathbb{R}$ if and only if $\varliminf\frac{b_n}{n}=0$ and $\varliminf\frac{c_n}{n}=0$. \qed 
\end{thm}

\section{Summary of all cases}
We are now in the position that we can describe all possible cases regarding $AAR_{(a_n)}$. Let us denote the accumulation points of $(a_n)$ by $(a_n)'$.

\begin{enumerate}

\item There is no finite accumulation point of $(a_n)$

\begin{enumerate}
\item If $(a_n)'=\{-\infty\}$ then $AAR_{(a_n)}=\{-\infty\}$.
\item If $(a_n)'=\{+\infty\}$ then $AAR_{(a_n)}=\{+\infty\}$.
\item If $(a_n)'=\{\pm\infty\}$ let $(a_n)=(b_n)||(c_n),\ b_n\to-\infty,\ c_n\to+\infty$. Then $AAR_{(a_n)}=\bar{\mathbb{R}}$ if and only if $\varliminf\frac{b_n}{n}=0$ and $\varliminf\frac{c_n}{n}=0$. Otherwise $AAR_{(a_n)}=\{\pm\infty\}$.
\end{enumerate}

\item If $\pm\infty\not\in (a_n)'$ then let $\varliminf a_n=a,\ \varlimsup a_n=b\ (a,b\in\mathbb{R})$. In this case $AAR_{(a_n)}=[a,b]$.

\item There are also finite and infinite accumulation points of $(a_n)$. Let $\varliminf a_n=\alpha\in\bar{\mathbb{R}},\ \varlimsup a_n=\beta\in\bar{\mathbb{R}},\ \min (a_n)'-\{-\infty\}=a\in\bar{\mathbb{R}},\ \max (a_n)'-\{+\infty\}=b\in\bar{\mathbb{R}}$. Obviously either $\alpha=-\infty$ or $\beta=+\infty$.

Then we can decompose $(a_n)$ in the following way: $(a_n)=(b_n)\cup(c_n)\cup(d_n),\ b_n\to\alpha,\ c_n\to\beta,\ (d_n)$ is the rest of the sequence (if needed). This decomposition is not unique in general. However if $\alpha=-\infty,\ a>-\infty$ then $(b_n)$ is unique up to finitely many elements and if $\beta=+\infty,\ b<+\infty$ then $(c_n)$ is unique up to finitely many elements.

\begin{enumerate}
\item If $a=-\infty$ or $b=+\infty$ then  $AAR_{(a_n)}=[a,b]\cup\{\alpha,\beta\}$.
\item If $-\infty<a\leq b<+\infty$ then we have the following cases:
\begin{enumerate}
\item Let $\alpha=-\infty,\ \beta<+\infty$. If $(b_n)$ is balanced then $AAR_{(a_n)}=[-\infty,b]$. If $(b_n)$ is not balanced then $AAR_{(a_n)}=[a,b]\cup\{-\infty\}$.
\item Let $\alpha>-\infty,\ \beta=+\infty$. If $(c_n)$ is balanced then $AAR_{(a_n)}=[a,+\infty]$. If $(c_n)$ is not balanced then $AAR_{(a_n)}=[a,b]\cup\{+\infty\}$.
\item Let $\alpha=-\infty,\ \beta=+\infty$. If $(b_n),(c_n)$ are both balanced then $AAR_{(a_n)}=\bar{\mathbb{R}}$. If $(b_n)$ is balanced, $(c_n)$ is not balanced then $AAR_{(a_n)}=[-\infty,b]\cup\{+\infty\}$. If $(b_n)$ is not balanced, $(c_n)$ is balanced then $AAR_{(a_n)}=[a,+\infty]\cup\{-\infty\}$. If none of $(b_n)$, $(c_n)$ is balanced then $AAR_{(a_n)}=[a,b]\cup\{-\infty,+\infty\}$.
\end{enumerate}
\end{enumerate}

\end{enumerate}

\section{On accumulation points}

\begin{thm}\label{tacc1}Let $(a_n)$ be a sequence bounded from below or above, $(a'_n)$ be one of its rearrangements. Let $c_n=\frac{\sum\limits_{i=1}^na'_i}{n}$. Then the set of accumulation points of $(c_n)$ is $[\varliminf c_n,\varlimsup c_n]$.
\end{thm}
\begin{proof}If $\varliminf c_n=\varlimsup c_n$ then there is nothing to prove. Assume that $\varliminf c_n<\varlimsup c_n$. Assume that $(a_n)$ is bounded from below and $a_n>K$. The other case is similar.

Let $\varliminf c_n<p<\varlimsup c_n$. Obviously it true for infinitely many $n$ that $c_n\leq p<c_{n-1}$. Then
\[\frac{a'_1+\dots+a'_n}{n}\leq p\]
which gives that
\[c_{n-1}=\frac{a'_1+\dots+a'_{n-1}}{n-1}\leq \frac{np-a'_n}{n-1}.\]
Then we get that
\[c_{n-1}-c_n=\frac{\sum\limits_{i=1}^{n-1}a'_i}{n-1}-\frac{\sum\limits_{i=1}^{n-1}a'_i+a'_n}{n}=
\frac{\sum\limits_{i=1}^{n-1}a'_i}{n(n-1)}-\frac{a'_n}{n}=\frac{c_{n-1}}{n}-\frac{a'_n}{n}\leq\]
\[\leq\frac{p}{n-1}-\frac{a'_n}{n(n-1)}-\frac{a'_n}{n}=\frac{p}{n-1}-\frac{a'_n}{n-1}\leq\frac{p-K}{n-1}<\epsilon\]
if $n$ is big enough.
\end{proof}

\begin{lem}\label{ll1}Let $(b_n)$ be a sequence, $c\in\mathbb{R}$. Assume that $\A(b_1,\dots,b_n)\in(a,b)$ if $n>N$ ($a<b,a,b\in\mathbb{R}$). Then $\forall\epsilon>0$ we can create a new sequence $(d_n)$ with $d_i=b_i\ (i<k),d_k=c,d_i=b_{i-1}\ (i>k>N)$ such that $n\geq k$ implies that $\A(d_1,\dots,d_n)\in(a-\epsilon,b+\epsilon)$.
\end{lem}
\begin{proof}Choose $k\in\mathbb{N}$ such that $n\geq k-1$ implies that 

(1) $|\frac{c}{n}|<\frac{\epsilon}{2}$ and

(2) $a-\frac{\epsilon}{2}<a\cdot\frac{n-1}{n}$.

If $m\geq k$ then \[\A(d_1,\dots,d_m)=\frac{c+\sum\limits_{i=1}^{m-1}b_i}{m}=\frac{c}{m}+\A(b_1,\dots,b_{m-1})\frac{m-1}{m}\] 
hence $a-\epsilon<\A(d_1,\dots,d_m)<b+\epsilon$.
\end{proof}

\begin{lem}\label{ll2}Let $(b_n)$ be a bounded sequence, $a=\varliminf b_n, b=\varlimsup b_n, a<b,\ (a,b\in\mathbb{R})$. Let $c,d\in\mathbb{R}$ such that $a<c<d<b$. Then there is a rearrangement $(b'_n)$ of $(b_n)$ and there is $N\in\mathbb{N}$ such that $n>N$ implies that $\A(b'_1,\dots,b'_n)\in(c,d)$.
\end{lem}
\begin{proof}It is an obvious consequence of \cite[Proposition 2.6]{lapaar}.
\end{proof}

\begin{thm}\label{taccmain}Let $(a_n)$ be a sequence such that it has at least 4 accumulation points: $a,b,-\infty,+\infty\ (a,b\in\mathbb{R},a<b)$. Let $Z\subset[a,b]$ be a closed set. Then there is a rearrangement $(a'_n)$ of $(a_n)$ such that the accumulation points of the sequence $(p_n)$ is exactly $Z$ where $p_n=\A(a'_1,\dots,a'_n)$.
\end{thm}
\begin{proof}Let $(a_n)=(b_n)||(c_n)||(d_n)||(e_n)||(f_n)$ where $b_n\to a,c_n\to b,d_n\to -\infty,e_n\to +\infty$ and $(f_n)$ is the rest of the elements in $(a_n)$. Assume that $\forall n\ a-1<b_n,c_n<b+1$. Let $K=\max\{|a-1|,|b+1|\}$.

First let us note that there is a sequence $(t_n)$ such that the accumulation points of $(t_n)$ is exactly $Z$. We know that $Z$ is separable because $\mathbb{R}$ is hereditary separable, hence let $Z_0\subset Z$ such that $cl(Z_0)=Z$ and $|Z_0|\leq\aleph_0$. Let $Z_0=\{w_n:n\in\mathbb{N}\}$ and then the sequence $(w_1,w_2,w_1,w_2,w_3,w_1,w_2,w_3,w_4,\dots)$ satisfies the requirement.

We are going to create a rearrangement $(a'_n)$ of $(a_n)$ such that there is a sequence $(n_k)$ such that $n_k\in\mathbb{N}, (n_k)$ is increasing and $n_k\leq n<n_{k+1}$ implies that $a_n\in(t_k-\frac{1}{k},t_k+\frac{1}{k})$. It is easy to see that $(a'_n)$ fulfills the requirements.

Let us generalize the problem slightly. We have a sequence of disjoint finite open intervals $(I_i)$ with $I_i\subset(a,b)$ and want to create a rearrangement $(a'_n)$ of $(a_n)$ such that there is a sequence $(n_k)$ such that $n_k\leq n<n_{k+1}$ implies that $a_n\in I_k$. 

We define sequences by recursion. By \ref{ll2} we can find a rearrangement $(r^1_n)$ of $(b_n)||(c_n)$ and $N_1\in\mathbb{N}$ such that $n\geq N_1$ implies that $\A(r^1_1,\dots,r^1_n)\in I_1$. We can also assume that $r^1_1=a_1$.

Assume that for $k\in\mathbb{N},k>1$ there are a sequence $(r^k_n)$ and positive integers $n_1<n_2<\dots<n_{k-1}<n_k$ such that 

(0) if $n>n_{k}$ then $r^k_n$ is an element of $(b_n)||(c_n)$,

(1) and if $n\leq n_{k-1}$ then $r^{k-1}_n=r^k_n$,

(2) and if $1\leq i\leq k-1$ and $n_i\leq n<n_{i+1}$ then $\A(r^i_1,\dots,r^i_n)\in I_i$, 

(3) and if $n_k\leq n$ then $\A(r^k_1,\dots,r^k_n)\in I_k$, 

(4) and there is $i\leq n_k$ such that $r^k_i=a_k$.

We will see that $n_1,\dots,n_k$ can be chosen independently of $k$.
\medskip

Assume that $\sup I_k<\inf I_{k+1}$. (The other inequality can be handled similarly.) 

If $a_{k+1}$ is already in $(r^k_n)$, say $r^k_i=a_{k+1}$, then let $N_k=\max\{i,n_k\}$. If not, then let $I_k=(a',b'),\ \epsilon=\frac{b'-a'}{4},\ I'=(a'+\epsilon,b-\epsilon)$ and apply \ref{ll1} 
for $(r^k_n),\ I',\ N=n_k,\  \epsilon,\ c=a_{k+1}$. Hence in that way we can merge $a_{k+1}$ into $(r^k_n)$ in a way that all conditions (0),(1),(2),(3),(4) remain valid (maybe we have to modify $n_k$). Let us keep the notation $(r^k_n)$ for this new sequence too. Let $N_k$ be the index of $a_{k+1}$ in $(r^k_n)$.

Let $h$ be the length of $I_{k+1}$. Let $M>\max\{N_k,\frac{9K}{h},\frac{18(b-a+1)}{h}\}$. If $m\geq M$ then $\frac{9(b-a+1)}{h}<\frac{m}{2}<\frac{m^2}{m+1}$ which gives that $(b-a+1)(m+1)<\frac{h}{9}m^2$. I.e. if we take intervals $\big((b-a+1)m,\frac{h}{9}m^2\big)$ for all such $m$, then those intervals will overlap. Hence there are a $P>M$ and an unused element $e$ from $(e_n)$ such that $(b-a+1)P<e<\frac{h}{9}P^2$ and $\frac{P-1}{P}a>a-1$.

First note that $b<\A(r^k_1,\dots,r^k_{P-1},e)$ because 
\[\A(r^k_1,\dots,r^k_P,e)=\frac{(P-1)\A(r^k_1,\dots,r^k_{P-1})+e}{P}>\frac{P-1}{P}a+\frac{e}{P}>b.\]
If $m>P$ then let $s=\sum\limits_{i=1}^{m-1}r^k_i,\ z=r^k_{m}$. Then
\[|\A(r^k_1,\dots,r^k_{m-1},e)-\A(r^k_1,\dots,r^k_{m},e)|=\left\lvert\frac{s+e}{m}-\frac{s+z+e}{m+1}\right\lvert=\left\lvert\frac{s+e-mz}{m(m+1)}\right\lvert\leq\]
\[\frac{m-1}{m(m+1)}|\A(r^k_1,\dots,r^k_{m-1})|+\frac{e}{m(m+1)}+\frac{|z|}{m+1}<\frac{K}{m}+\frac{e}{P^2}+\frac{K}{m}<\frac{3h}{9}=\frac{h}{3}.\]
Now there is a sufficiently large $m$ such that $m>P$ and 
\[\A(r^k_1,\dots,r^k_P,e,r^k_{P+1},\dots,r^k_m)\in I_{k+1}\] because $\A(r^k_1,\dots,r^k_P,e)>b$ and
\[\varlimsup\limits_{m\to\infty}\A(r^k_1,\dots,r^k_P,e,r^k_{P+1},\dots,r^k_m)\leq\sup I_k<\inf I_{k+1}\] and by $|\A(r^k_1,\dots,r^k_{m-1},e)-\A(r^k_1,\dots,r^k_{m},e)|<\frac{h}{3}$. Fix such $m$.

Let 
\[r^{k+1}_n=
\begin{cases}
r^k_n&\text{if }n\leq m\\
e&\text{if }n=m
\end{cases}.\]
If $n>m$ then choose not-used elements from $(b_n)||(c_n)$ such that \[\A(r^{k+1}_1,\dots,r^{k+1}_n)\in I_{k+1}.\] We show that this can be done (we apply a similar argument as above). We go by induction. It is true for $n=m$. Assume that we are done till $n$ and looking for $r^{k+1}_{n+1}$. Let  $s=\sum\limits_{i=1}^{n-1}r^{k+1}_i,\ z\in(a-1,b+1)$.
\[|\A(r^{k+1}_1,\dots,r^{k+1}_{n})-\A(r^{k+1}_1,\dots,r^{k+1}_{n},z)|=\left\lvert\frac{s}{n}-\frac{s+z}{n+1}\right\lvert=\left\lvert\frac{s-nz}{n(n+1)}\right\lvert\leq\]
\[\frac{1}{n+1}|\A(r^{k+1}_1,\dots,r^{k+1}_{n})|+\frac{|z|}{n+1}<\frac{K}{n}+\frac{K}{n}<\frac{2h}{9}=\frac{h}{2}.\]
If $\A(r^{k+1}_1,\dots,r^{k+1}_{n})\in I_{k+1}$, say $\A(r^{k+1}_1,\dots,r^{k+1}_{n})\leq\inf I_{k+1}+\frac{h}{2}$, then we can choose $r^{k+1}_{n+1}$ from $(c_n)$ such that $r^{k+1}_{n+1}>\A(r^{k+1}_1,\dots,r^{k+1}_{n})$ and that will satisfy $A(r^{k+1}_1,\dots,r^{k+1}_n)\in I_{k+1}$. The other case is similar.

\medskip

Let now $n_{k+1}=m$. We show that all conditions (0),(1),(2),(3),(4) will be valid for $(r^{k+1}_n)$. (0) and (1) are obvious by definition. (2) is valid because we keep all previous $n_i\ (1\leq i\leq k)$ and by (3) for $(r^{k}_n)$. (3) holds by definition. The first step was to include $a_{k+1}$ hence (4) is valid too.

Let $a'_n=r^k_n$ if $n_k\leq n<n_{k+1}$. Clearly the accumulation points of $p_n=\A(a'_1,\dots,a'_n)$ is exactly $Z$. 
\end{proof}

\begin{rem}In \ref{taccmain} if we want to weaken the condition then by \ref{tacc1} we cannot abandon neither $-\infty$ nor $+\infty.$
\end{rem}

\begin{rem}In \ref{taccmain} we cannot say more in general in a way that any closed sets $Z$ outside $[a,b]$ could be the accumulation points of the averages of a rearranged sequence. 
\end{rem}
\begin{proof}
Let $(a_n)=(-2^n)||(-1)||(1)||(2^n)$. By \ref{t1aariff2b} and \ref{t2t1} $AAR_{(a_n)}=[-1,1]\cup\{-\infty,+\infty\}$.
\end{proof}


{\footnotesize

\noindent

\noindent E-mail: alosonczi1@gmail.com\\
}

\begin{thebibliography}{www}


\bibitem{be} F. Bagemihl, P. Erd\H os, {\em Rearrangements of $C_1$-summable series}, Acta Math. {\em 92} (1954), 35--53.

\bibitem{bullen} P. S. Bullen, {\em Handbook of means and their inequalities}, vol. 260 Kluwer Academic Publisher, Dordrecht, The Netherlands (2003).

\bibitem{ds} Y. Dybskiy, K. Slutsky, {\em Riemann Rearrangement Theorem for some types of convergence}, preprint, arXiv:math/0612840

\bibitem{fs} R. Filipów, P. Szuca, {\em Rearrangement of conditionally convergent series on a small set}, J. Math. Analysis and Appl. {\em 362 (1)} (2010), 64--71.

\bibitem{lz} G. G. Lorentz, K. I. Zeller, {\em Series rearrangements and analytic sets}, Acta Math. {\em 100} (1958), 149--169.


\bibitem{lapaar} A. Losonczi, {\em Points accessible in average by rearrangement of sequences I}, Transnational Journal of Mathematical Analysis and Applications Volume 7, Issue 1 (2019), 1 -- 27

\bibitem{p} P. A. B. Pleasants, {\em Rearrangements that preserve convergence}, J. London Math. Soc.  {\em 15} (1977),  134--142.

\bibitem{s} M. A. Sarigöl, {\em Rearrangements of bounded variation sequences}, Proc. Indian Acad. Sci. (Math. Sci.) {\em 104 (2)} (1994), 373--376.

\end{thebibliography}
\end{document}